\documentclass[a4paper,12pt]{article} 
\newdimen\paperhight
\usepackage{amsmath}
\usepackage{amssymb}
\usepackage{amsfonts}
\usepackage[usenames]{color}

\newcommand{\dsp}{\displaystyle}
\newcommand{\dto}{\downarrow}

\newcommand{\pr}{\par \vspace{3mm}\noindent [{\bf Proof}] \qquad}
\newcommand{\prend}{\hfill \qed \par \vspace{3mm}}

\newcommand{\qed}{\quad\hbox{\rule[-2pt]{3pt}{6pt}}\par\vspace{3mm}}

\newcommand{\1}{{\bf 1}} 
 
\newcommand{\C}{\mathbb C} 
\newcommand{\Z}{\mathbb Z} 
\newcommand{\Q}{\mathbb Q} 
 
\newcommand{\N}{\mathbb N}

\newcommand{\CI}{{\cal I}}
\newcommand{\CF}{{\cal F}}

\newcommand{\CU}{{\cal U}}

\newcommand{\CY}{{\cal Y}}

\newcommand{\al}{\alpha}
\newcommand{\be}{\beta}

\newcommand{\ga}{\gamma}
\newcommand{\la}{\lambda}

\newcommand{\wt}{{\rm wt}}
\newcommand{\Hom}{{\rm Hom}}

\newcommand{\End}{{\rm End}}

\newcommand{\Ker}{{\rm Ker}}
\newcommand{\Res}{{\rm Res}}

\newtheorem{thm}{Theorem}
\newtheorem{prn}[thm]{Proposition}
\newtheorem{dfn}[thm]{Definition}

\newtheorem{lmm}[thm]{Lemma}
\newtheorem{rmk}[thm]{Remark}

\begin{document}
\title{Flatness of Tensor Products and 
Semi-Rigidity for \\
$C_2$-cofinite Vertex Operator Algebras I}
\author{\begin{tabular}{c}
Masahiko Miyamoto\\
Department of Mathematics, \\
University of Tsukuba, \\
Tsukuba, 305 Japan \end{tabular}}
\date{}
\maketitle

\abstract{We study properties of a $C_2$-cofinite vertex 
operator algebra $V=\oplus_{i=0}^{\infty}V_i$ of 
CFT type.  If it is also rational (i.e. all modules are completely reducible) and 
$V'\cong V$, then the rigidity of the tensor category 
of modules has been proved by Huang \cite{H2}, where $V'$ denotes the restricted dual 
of $V$.  However, when we treat 
irrational $C_2$-cofinite VOAs, the rigidity is too strong, 
because it is almost equivalent to be rational as we see.  
We introduce a weaker condition "semi-rigidity". 
We expect that all $C_2$-cofinite VOAs satisfy this condition. 
Under the assumption of the semi-rigidity and the existence of canonical homomorphisms, 
we prove the following results. 
We show that if $P$ is a projective cover of a $V$-module $V$, then 
for any finitely generated $V$-module $M$, 
its projective cover is a direct summand of 
the tensor product $P\boxtimes M$ (defined by logarithmic intertwining operators) 
of $M$ and $P$. Using this result, we prove the flatness property of 
finitely generated modules for the tensor products, that is, 
if $0\rightarrow A\rightarrow B\rightarrow C\rightarrow 0$ is exact then so is $0\rightarrow D\boxtimes A\rightarrow D\boxtimes B\rightarrow 
D\boxtimes C\rightarrow 0$ for any finitely generated $V$-modules $A$, $B$, $C$ and $D$.
As a corollary, we have that 
if a semi-rigid $C_2$-cofinite $V$ contains a rational subVOA with the same Virasoro element, 
then $V$ is rational.}

\section{Introduction}
The concept of a vertex (operator) algebra 
(shortly VOA) $V=(V,Y,{\bf 1},\omega)$ 
was introduced 
by Borcherds \cite{B} and then by Frenkel-Lepowsky-Meurman \cite{FLM} 
with a motivation to explain a mysterious relation between the monster simple group (the largest 
sporadic finite simple group) and 
the $j$-function ($j(z)=q^{-1}+744+196884q+\cdots$). 
A vertex algebra was introduced by axiomatizing the concept of a 
Chiral algebra and so it offers a rigorous proof in the two-dimensional Conformal Field Theory. 
Furthermore, because it has infinitely many operators, 
it has rich connection with algebraic objects. 

Classically, because of the physical understanding and technical reason, 
the objects interested in are mainly $\N$-graded simple modules or 
vertex operator algebras whose modules are completely reducible.  
Such vertex operator algebras are called rational. 
However, the research for irrational 
vertex operator algebras has recently been progressive and 
ring theoretical approaches are getting more important. 
For example, a Zhu algebra $A_0(V)$ introduced by 
Zhu \cite{Z} and $n$-th Zhu algebras $A_n(V)$ 
extended by Dong-Li-Mason \cite{DLM} 
offer a lot of information of modules.

Among irrational VOAs, if $V$ satisfies 
$C_2$-cofiniteness condition introduced by Zhu \cite{Z}, 
then $V$ has only finitely many simple modules and 
satisfies nice properties (see Proposition \ref{C2Pr}) including a modular 
invariance property of (pseudo) trace functions \cite{Miy}.  
Here, a $V$-module $W$ is called $C_2$-cofinite if $\dim W/C_2(W)<\infty$, 
where $C_2(W)=\langle v_{-2}w \mid v\in V, \wt(v)> 0, w\in W\rangle$. 
Daring to say too much, under the $C_2$-cofiniteness condition, 
a vertex operator algebra seems to satisfy the most of 
nice properties which hold in the standard theories  
on groups and finite dimensional rings. 
For example, as we will explain, 
for two $V$-modules $W$ and $U$, its tensor product 
$W\boxtimes U$ (defined by logarithmic intertwining operators) is 
well-defined. 
In this paper, we will not use any explicit construction of 
tensor products, but a property that it is a largest one in a sense. 
Although Huang, Lepowsky and Zhang have extended their $P(z)$-tensor product theory to 
logarithmic intertwining operators for more general setting in \cite{HLZ}, 
as long as we consider a $C_2$-cofinite VOA, we will choose an easier way to treat them. 
Let us explain tensor products from this point of view. 
For any f.g. modules $W$ and $U$, 
we can consider the set $\CF(W,U)$ of pairs $(T,\CY)$ of a f.g. module $T$ and 
a surjective (logarithmic) intertwining operator $\CY\in \CI_{W,U}^{T}$. 
In this case, a tensor product $(W\boxtimes U,\CY^{W\boxtimes U})$ of 
two modules $W$ and $U$ should be understood as 
an isomorphic class satisfying the universal property, in other words, 
for any element $(T,\CY)\in \CF(W,U)$ there is a homomorphism 
$\phi\in \Hom(W\boxtimes U, T)$ such that 
$\phi\cdot\CY^{W\boxtimes U}=\CY$.  
As we will explain later, $\CF(W,U)$ becomes a (right) direct set, 
that is, it has a partial order $<$ such that for any two elements $\al$ and $\be$ there 
is a larger element $\gamma$ satisfying $\al<\ga$ and $\be<\ga$. 
Since $\CF(W,U)$ contains only finitely many non-isomorphic classes 
when $V$ is $C_2$-cofinite as we will see, 
$\CF(W,U)$ contains a unique maximal element up to isomorphism. 
This is a "tensor product" of $W$ and $U$. From the maximality of tensor products, 
for any homomorphism $\tau:A\rightarrow B$, we can induce a canonical homomorphism 
$\tau\boxtimes {\rm id}_D:A\boxtimes D\rightarrow B\boxtimes D$, where 
${\rm id}_D$ denotes the identity on $D$, (see the statement before Proposition \ref{RightFl}). 
For the tensor products, there are canonical isomorphisms 
$$\begin{array}{c}
\sigma_{W,U}:W\boxtimes U \rightarrow U\boxtimes W \cr
\mu:(W\boxtimes U) \boxtimes T \rightarrow W\boxtimes (U\boxtimes T),
\end{array}$$
see \cite{H1}.  What we need is the existence of such canonical homomorphism and so 
we will use these results without giving explanation.

Throughout this paper, we assume the existence of dual elements in the category of $V$-module. 
Namely, we assume that for any $V$-module $W$ there are epimorphism 
$$e_W:\widetilde{W}\boxtimes W\rightarrow V \mbox{  and  } e_{\widetilde{W}}:W\boxtimes \widetilde{W}\rightarrow V, 
\eqno{(1.1)}$$
for some $V$-module $\widetilde{W}$, where $e_{\widetilde{W}}=e_W\sigma_{W,\widetilde{W}}$.  
For example, if $V\cong V'$, then we can take $\widetilde{W}=W'$, where $W'$ denotes 
a restricted dual of $W$, (see (2.1) or \cite{FHL} for its definition.)  
On the other hand, as we explained, there is a canonical isomorphism 
$$\mu:(W\boxtimes \widetilde{W})\boxtimes W  \rightarrow W\boxtimes (\widetilde{W} \boxtimes W).
 \eqno{(1.2)}$$
The rigidity of tensor category is defined by 
$$({\rm id}_W\boxtimes e_W)\cdot\mu\cdot(i_W\boxtimes {\rm id}_W)={\rm id}_W$$ 
under the assumption of the existence of $i_W:V \rightarrow W\boxtimes \widetilde{W}$,   
where we identify $W\boxtimes V$ and $V\boxtimes W$ with $W$.  
If $V$ is rational, we know from the definition that the rigidity of $V$ is 
equivalent to nonvanishing of 
the denominators of Verlinde formulas. The latter statement was proved by 
Huang in the proof of Verlinde formula if $V$ is rational, 
$C_2$-cofinite and $V\cong V'$ (see \cite{H2}).  

Because our target contains non-semisimple modules, 
$e_{\widetilde{W}}:W\boxtimes \widetilde{W}\rightarrow V$ may not be split and so 
we can't expect to have an embedding $i_W:V\to W\boxtimes \widetilde{W}$. Therefore, 
we will consider a homomorphism $\rho:P\to W\boxtimes \widetilde{W}$ such that 
$e_{\widetilde{W}}(\rho(P))=V$, where $P$ is a projective cover of $V$. 
Now we introduce the following weaker condition:

\begin{dfn}
We will call that a f.g. $V$-module $W$ satisfies {\bf semi-rigidity} if there exist 
the following diagram 
$$ \begin{array}{rrlcrl}
&P&\boxtimes W && &\cr
&\rho \downarrow & \boxtimes \downarrow id_W&  & &\cr
\mu:&(W\boxtimes \widetilde{W})&\boxtimes W&\xrightarrow{\cong}&W\boxtimes&(\widetilde{W}\boxtimes W) \cr 
&e_{\widetilde{W}}\downarrow &\boxtimes \downarrow id_W&&id_W \downarrow \boxtimes &\downarrow e_W\cr
&V&\boxtimes W&& W\boxtimes &V 
\end{array}$$
such that $e_{\widetilde{W}}\rho(P)=V$ and 
$$({\rm id}_W\boxtimes e_W)(\mu(\rho(P)\boxtimes W))\cong W, $$ 
where $\widetilde{W}$ is a f.g.$V$-module, 
$e_{\widetilde{W}}:W\boxtimes \widetilde{W}\rightarrow V$ 
and $e_W:\widetilde{W}\boxtimes W\rightarrow V$ are given epimorphisms, and 
$P$ is a projective cover of $V$. 
In particular, $V$ is called semi-rigid if all f.g. $V$-modules satisfy the semi-rigidity.
\end{dfn}

The main purpose in this paper is to prove that the flatness of modules for 
tensor products under the assumption of the semi-rigidity. 
Namely, we will prove the following theorem in \S 3.3.\\ 

\noindent
{\bf Theorem \ref{MainTh}}\quad {\it Let $V$ be a $C_2$-cofinite 
vertex operator algebra of 
CFT type and assume that all simple $V$-modules satisfy the semi-rigidity. 
If $A, B, C$, and $D$ are f.g. $V$-modules and 
$$0 \rightarrow A \xrightarrow{\tau} B \xrightarrow{\sigma} C\rightarrow 0$$ 
is an exact sequence of $V$-modules, then so is 
$$0\rightarrow A\boxtimes D \xrightarrow{\tau\boxtimes {\rm id}_D} 
B\boxtimes D 
\xrightarrow{\sigma\boxtimes {\rm id}_D} C\boxtimes D \rightarrow 0. $$}

As we will show later, under $C_2$-cofiniteness condition, 
every f.g. module has a projective cover. In particular, we call 
a projective cover $\rho:P\rightarrow V$ (or a module $P$ itself) 
of a $V$-module $V$ "a principal projective cover". 
The key result in the proof of the above theorem is the following, 
which will be proved in \S 3.1. \\

\noindent
{\bf Theorem \ref{KeyTh}}\quad {\it  
Assume that $V$ is a $C_2$-cofinite VOA of CFT type and all simple $V$-modules satisfy 
the semi-rigidity.  
Let $W$ be a f.g. $V$-module and 
$f:S \rightarrow W$ a projective cover of $W$.  
Then there is an epimorphism $g:W\boxtimes P \rightarrow S$ such that 
$fg$ coincides with ${\rm id}_W\boxtimes \rho:W\boxtimes P \rightarrow W\boxtimes V=W$, 
where $P$ denotes a principal projective cover. 
where $P$ denotes a principal projective cover. 
In particular, $S$ is a direct summand of $W\boxtimes P$. }\\

As easy corollaries of these theorems, we will have the followings:\\

\noindent
{\bf Corollary 16}\quad {\it
Let $V$ be a $C_2$-cofinite vertex operator algebra of CFT type and 
assume that $V$ satisfies the semi-rigidity. 
If $V$ is projective as a $V$-module, then $V$ is rational. }\\

\noindent
{\bf Corollary 17}\quad {\it
Let $V$ be a $C_2$-cofinite vertex operator algebra of CFT type and 
assume that $V$ satisfies the semi-rigidity. 
If $V$ contains a rational subVOA $W$ with the same Virasoro element, then 
$V$ is rational.  }\\

It is easy to see that if a $C_2$-cofinite subVOA $U$ of $V$ satisfies the semi-rigidity 
then $V$ satisfies the semi-rigidity. Therefore, we have: \\

\noindent
{\bf Corollary 18}\quad {\it 
Let $V$ be a vertex operator algebra of CFT type and 
$U$ a subVOA of $V$ with the same Virasoro element. 
If $U$ is $C_2$-cofinite, rational and rigid, 
then so is $V$. }\\

\noindent
{\rm Acknowledgement} \\
The author wishes to thank S.~Naito and D.~Sagaki 
for their advice. 
The author also would like to express 
the special appreciation to T.~Abe and A.~Hida for their help. \\

\section{Preliminary results}
\subsection{Notation}
  Throughout this paper, all vector spaces are over the complex number field $\C$ and 
$V$ denotes a vertex operator algebra $(V,Y(\cdot,z),\1,\omega)$, 
where $V=\oplus_{n=0}^{\infty}V_n$, $\omega$ is a Virasoro element, 
$\1$ is the vacuum, and $Y(v,z)=\sum_{n\in \Z} v_nz^{-n-1}\in 
\End(V)[[z,z^{-1}]]$ denotes a 
vertex operator of $v\in V$. They satisfy the conditions 
$1\sim 4$ in Definition 2 by replacing all $W,U,T$ by $V$. We also have 
$Y(\1,z)={\rm id}_V$ and the coefficients of $Y(\omega,z)$ satisfy the Virasoro 
algebra relations. If $\dim V_0=1$, then $V$ is called CFT-type. 
If $W=V$, $U=T=M$ in Definition 2 and a set  
$$\{Y^M(v,z)=\sum_{n\in \Z} v^M_nz^{-n-1}\in 
\End(M)[[z,z^{-1}]] \mid v\in V\}$$ 
satisfies the same conditions $1\sim 4$, then $M$ is called 
a (weak) $V$-module. We note that if $V$ is $C_2$-cofinite, then $V$ has 
only finitely many isomorphism classes of simple $V$-modules and 
all weak modules are $\N$-gradable (see Proposition \ref{C2Pr}). Hence 
a module $M$ is finitely generated if and only if $M$ has a composition 
series of finite length.  
We also note that for the degree operator $L(0)=\omega^M_1$ 
on a $V$-module $M$, 
a homogeneous space $M_r$ with eigenvalue $r$ is not 
necessary to be an eigenspace, but a generalized eigenspace of $L(0)$. Namely, for $r\in \C$  
$$M_r=\{w\in M \mid (L(0)-r)^Kw=0 \mbox{ for some }K\}.$$ 
We use notation weight $"\wt"$ to denote eigenvalues (or the semisimple part) of $L(0)$, 
that is, $\wt(d)=r$ for $d\in M_r$. 
The lowest weight of a module $M$ is called a conformal weight of $M$ and 
denoted by $\wt(M)$. It is known that the weights of all modules 
are rational numbers if $V$ is $C_2$-cofinite \cite{Miy}. 
For a $\Q$-graded module $M=\oplus_{n\in \Q}M_r$, 
$M'$ denotes the restricted dual $V$-module 
$\oplus_{r\in \Q}(\Hom(M_r,\C))$, where an adjoint vertex operator 
$Y^{M'}(v,z)$ on $M'$ is given by 
$$\langle Y^{M'}(v,z)w',w\rangle
=\langle w',Y(e^{zL(1)}(-z^{-2})^{L(0)}v,z^{-1})w\rangle 
\eqno{(2.1)}$$ 
for $w'\in M'$ 
and $w\in M$ and $\langle w',w\rangle$ denotes $w'(w)\in \C$, 
see \cite{FHL}.

\subsection{(logarithmic) intertwining operators}
Similar to a vertex operator 
$Y^M(v,z)=\sum v^M_nz^{-n-1}\in \End(M)[[z,z^{-1}]]$ for some module $M$, 
it is natural to consider an intertwining operator from 
a module $U$ to another module $T$ as a formal power series (\cite{TK}). 
However, without the assumption of rationality, 
there is no reason for an intertwining operator to have a specific 
form like a formal power series. 
Fortunately, if $V$ satisfies the $C_2$-cofiniteness 
condition, then all f.g. modules are 
$C_1$-cofinite. Hence, as Huang has shown \cite{H1}, 
for any f.g. modules $U,W,T$ and 
any intertwining operators $\CY(\ast,z)\in \CI_{W,U}^{T}$, 
its correlation function 
$$  f(u,w,t;\la, z)=\la^{\wt(w)+\wt(v)+\wt(u)-\wt(t)}\langle t', 
Y^T(v,\la)\CY(w,\la z)u\rangle$$
satisfies a differential equation of 
regular singular points for any $\la\in \C$, $u\in U$, $w\in W$, $t'\in T'$, $v\in V$. 
Hence we may assume that each $\CY(\ast,z)$ has a shape of formal power series with $\log z$ terms: 
$$\CY(w,z)=\sum_{i=0}^K\sum_{m\in {\mathbb C}} 
w^{\CY}_{(m,i)}z^{-m-1}\log^i z\in \Hom(U,T)\{z\}[\log z] $$
for all $w\in W$. 
Such a formal power series with (bounded) natural integer powers of 
$\log z$ is called "logarithmic type" (see \cite{Mil}). 
As long as we study an irrational $C_2$-cofinite VOA, 
it is natural to treat such all intertwining operators. 
Therefore, we will call them (including the case $K=0$, that is, 
ordinary intertwining operators of formal power series) intertwining operators 
in this paper.

In 1993, Gurarie \cite{G} construct a CFT-like structure with logarithmic behavior. 
Let us review the definition of logarithmic intertwining operators, 
see \cite{Mil} (and \cite{F}, too). 

\noindent
\begin{dfn} Let $W$, $U$ and $T$ be 
f.g. $\N$-graded $V$-modules. 
A (logarithmic) intertwining operator of type 
$(W:U\rightarrow T)$ is a linear map 
$$\begin{array}{l}
  \CY(,z):W \rightarrow \Hom(U,T)\{z\}[\log(z)] \cr
  \displaystyle{\CY(w,z)=\sum_{i=0}^k\sum_{r\in \C} 
w_{(r,i)}z^{-r-1}\log^iz} 
\end{array}$$
satisfying the following conditions: \\
1. The lower truncation property: for each $u\in U$ and $i$,  
$w_{(r,i)}u=0$  for ${\rm Re}(r)>\!\!>0$, 
where ${\rm Re}(r)$ denotes the real part of $r\in \C$.  \\
 2. $L(-1)$-derivative property: 
$\CY(L(-1)w,z)=\frac{d}{dz}\CY(w,z)$ for $w\in W$. \\
3. Commutativity: $v_n^T\CY(w,z)-\CY(w,z)v_n^U=
\displaystyle{\sum_{i=0}^{\infty}} \binom{n}{i}\CY(v_i^Ww,z)z^{n-i} 
\mbox{  for  }v\in V$. \\
4. Associativity: for $v\in V$ and $n\in \N$, \\
$\CY(v^W_{n}w,z)=
\displaystyle{\sum_{m=0}^{\infty}}(-1)^m\binom{n}{m}v^T_{-m-1}z^m\CY(w,z)
+\CY(w,z)
\displaystyle{\sum_{m=0}^{\infty}}(-1)^{m+n+1}\binom{n}{m}v^U_mz^{-m-1},$ 
\end{dfn}

\noindent
where $Y^X(v,z)=\sum v^X_nz^{-n-1}$ denotes a vertex operator of $v\in V$ on 
a $V$-module $X$. 
If $\CY$ does not contain $\log z$ terms, we call it "ordinary". 
The set of intertwining operators of type 
$(W:U\rightarrow T)$ becomes a vector space, which we denote by 
$\CI_{W,U}^T$. 
It is known that the above Commutativity and 
Associativity are replaced by Borcherds' identity: \\
For $v\in V$, $w\in W$, $u\in U$ and $p,q\in\Z$ and $n\in\C$, we have
$$ \begin{array}{l}
\displaystyle{\sum_{i=0}^{\infty}\binom{q}{i}(v^W_{(p+i)}w)_{(q+n-i,r)}u}\cr
\hspace{2cm}=\displaystyle{\sum_{i=0}^{\infty}(-1)^i
\binom{p}{i}(v^T_{(p+q-i)}w_{(n+i,r)}u-(-1)^pw_{(p+n-i,r)}v^U_{(q+i)}u)}.
\end{array}\eqno{(2.2)}$$

Let $\CY^{(i)}(w,z)$ denote the coefficient term 
$\sum_{n\in \C} w_{(n,i)}z^{-n-1}$ of 
$\log^iz$ in $\CY(w,z)$ for $i=0,1,...,K$. 
We note that if $V$ is finitely generated, 
then there is an integer $K$ such that 
$\CY^{(n)}(w,z)=0$ for any $n>K, w\in W$ and $\CY$. 
In this paper, $K$ denotes the largest integer such that 
$\CY^{(K)}\not=0$. 
Moreover, since vertex operators $Y^M(v,z)$ on modules $M$ 
have no $"{\log} z"$ terms, 
every $\CY^{(i)}$ satisfies all properties of intertwining 
operators except the $L(-1)$-derivative property. 
On the other hand, from the $L(-1)$-derivative property for $\CY$, 
we have two important properties:
$$\begin{array}{rl}
\displaystyle{\CY^{(m)}(w,z)=}&\displaystyle{\frac{1}{m!}(z\frac{d}{dz}-zL(-1))^m\CY^{(0)}(w,z)}, 
\qquad \mbox{ and }\cr
\displaystyle{(i+1)w_{(n,i+1)}u
=}&\displaystyle{-(L(0)\!-\!\wt)w_{(n,i)}u+((L(0)\!-\!\wt)w)_{(n,i)}u+w_{(n,i)}((L(0)\!-\!\wt)u)}. 
\end{array} \eqno{(2.3)}$$
In particular, ${\CY}^{(K)}(\ast,z)$ is an ordinary 
intertwining operator (i.e. of formal power series).  
On the other hand, from (2.3), 
$$(z\frac{d}{dz}-zL(-1))^{K+1}\CY^{(0)}(w,z)=0 \eqno{(2.4)}$$
holds for all $w\in W$. 
Such a formal power series, that is, a logarithmic formal power series satisfying 
all conditions in Definition 2 except $L(-1)$-derivative property but (2.4) 
is called an $L(-1)$-nilpotent intertwining operator. 
Conversely, from such an $L(-1)$-nilpotent intertwining operator $\widetilde{\CY}_0$, 
we can construct a logarithmic intertwining operator 
$$\widetilde{\CY}(w,z)=\sum_{i=0}^{\infty} \left\{\frac{1}{i!}(zL(-1)-z\frac{d}{dz})^i\widetilde{\CY}_0(w,z)\right\}\log^iz.  \eqno{(2.5)}$$
The following comes from (2.3) easily. 

\begin{lmm} If $L(0)$ acts on $W,U,T$ semi-simply, then 
$\CY$ is an ordinary intertwining operator for $\CY\in \CI^T_{W,U}$. 
\end{lmm}

\subsection{$C_2$-cofiniteness and tensor products}
In this subsection, we explain about 
tensor products of modules defined by (logarithmic) 
intertwining operators under the assumption of $C_2$-cofiniteness.  

For a $V$-module $W$ and $m\in \N$, set 
$$C_m(W)=\langle v_{-m}w\mid v\in V, \wt(v)>0, w\in W\rangle. \eqno{(2.6)}$$
If $W/C_m(W)$ has a finite dimension, then we call $W$ to be $C_m$-cofinite. 
Among them, $C_2$-cofiniteness is the most important and
offers many nice properties. 
For example, we have:

\begin{prn}\label{C2Pr}
Let $V$ be a $C_2$-cofinite VOA. Then we have the followings: \\
(i) Every weak module is $\Z_+$-gradable and weights are all rational numbers. 
Moreover, this condition is equivalent to $C_2$-cofiniteness, {\rm \cite{Miy}}. \\
(ii) Evey $n$-th Zhu algebra $A_n(V)$ is finite dimension and
the number of inequivalent simple modules is finite, 
{\rm \cite{GN},\cite{DLM}}.\\
(iii) Set $V=B_2(V)+C_2(V)$ for 
a finite dimensional subspace $B_2(V)$ spanned by homogeneous elements. 
Then for any weak module $W$ generated from one element $w$ has the 
following spanning set $\{v^1_{n_1}....v^k_{n_k}w 
\mid v^i\in B_2(V), 
\quad  n_1<\cdots <n_k\}$.   
In particular, $\dim W_n$ is bounded by a function on $n$ and $\wt(w)$. 
Moreover, every f.g. $V$-module is $C_n$-cofinite for 
any $n=1,2,...$, {\rm \cite{Miy},\cite{Bu},\cite{GN}}. 
\end{prn}

For f.g. $V$-modules $W$ and $U$, consider 
the set of surjective intertwining operators  
$$ \CF(W,U)=\{(\CY,T)\mid  T 
\mbox{ is a f.g. $V$-module, $\CY\in \CI_{W,U}^T$ is "surjective"}\}, \eqno{(2.7)}$$
where $\CY(w,z)=\sum_{i=0}^K\sum_{m\in \C} w_{(m,i)}z^{-m-1}\log^iz$ is called surjective if 
$$<w_{(m,i)}u\mid w\in W, u\in U, m\in \C, i=0,...,K>=T.$$
For two surjective operators 
$\CY^1\in \CI_{W,U}^{T^1}$, $\CY^2\in \CI_{W,U}^{T^2}$,  
we define a partial order 
$\CY^1< \CY^2$ in $\CF(W,U)$ if there is a $V$-homomorphism 
$f\in \Hom(T^2,T^1)$ such that 
$$  f(\CY^2(w,z)u)=\CY^1(w,z)u  \quad {}^{\forall} w\in W, {}^{\forall} u\in U .$$
Clearly, if $\CY^1< \CY^2$ and $\CY^2< \CY^1$, then we 
have $T^1\cong T^2$ and we call them isomorphic. 

Let us show that $\CF(W,U)$ is a (right) directed set. 
For any $\CY^1,\CY^2\in \CF(W,U)$, say 
$\CY^1\in \CI_{W,U}^{T^1}$ and $\CY^2\in \CI_{W,U}^{T^2}$, 
we define $\CY$ by 
$$\CY(w,z)u=(\CY^1(w,z)u, \CY^2(w,z)u)\in (T^1\oplus T^2)\{z\}[\log z].$$
Clearly, $\CY$ is an intertwining 
operator. Let $T$ denote the subspace spanned by all images of $\CY$, 
then since $T_1$ and $T_2$ have composition series of finite length, 
so does $T$ and so $(\CY,T)\in \CF(W,U)$. 
By the projections $\pi_i:T\subseteq T_1\oplus T_2 \rightarrow T_i$, we have 
$\pi_i(\CY)=\CY^i$ for $i=1,2$. Namely, $\CY^1< \CY$, $\CY^2< \CY$.

\begin{dfn} 
If there is a unique maximal element $\CY\in \CI_{W,U}^T$ in 
$\CF(W,U)$ up to isomorphism, 
we call $T$ a tensor product of $W$ and $U$ and we denote it by $W\boxtimes U$. 
\end{dfn}

Since $\CF(W,U)$ has a larger element for any finite subsets, 
a necessary and sufficient 
condition for the existence of a f.g. tensor product module 
is that $\CF(W,U)$ contains only finitely many 
isomorphism classes.  \\

Before we start the proof of finiteness of $\CF(W,U)$, 
we give a brief review of an $n$-th 
Zhu algebra $A_n(V)$ from a view point of operators. 
For an $\N$-graded module $W=\oplus_{i=0}^{\infty} W(i)$, 
we restrict the grade preserving operators $o(v)=v^W_{\wt(v)-1}$ of 
$v\in V$ to the actions on $\oplus_{i=0}^nW(i)$. 
From the Borcherds' identity, for any $v,u\in V$, we can define 
an element $v\ast u\in V$ which does not depend on the 
choice of $W$ such that $o(v)o(u)=o(v\ast u)$ on $\oplus_{i=0}^nW(i)$. 
Therefore, if we set 
$$O_n(V)=\{v\in V\mid o(v)=0 \mbox{ on $\oplus_{i=0}^nW(i)$ for any 
$\N$-graded modules $W$}\},$$
then $A_n(V)=V/O_n(V)$ becomes an associative algebra. 
This is an $n$-th Zhu algebra. 
From Borcherds' identity (2.2), for any $v,u\in V$, we know 
$$ \sum_{i=0}^{\infty}\binom{\wt(v)+n}{i}v_{-2n-2+i}u \in O_n(V).\eqno{(2.8)}$$
The wonderful points done by 
Zhu \cite{Z} and Dong-Li-Mason \cite{DLM} are that 
they have determined $O_n(V)$ explicitly without using modules 
and then reconstructed $V$-modules from $A_n(V)$-modules. 

\begin{rmk}
Because we will treat reducible modules, we have to 
be careful about $n$-th homogeneous subspaces. If $U=\oplus_{i=0}^{\infty}U(i)$ is a 
submodule of $W=\oplus_{i=0}^{\infty}W(i)$, 
then an $n$-th homogeneous space $U(n)$ of $U$ and $U\cap W(n)$ may be different. 
In order to avoid this situation, we will consider fixed weights as follows:

From Proposition 5, $V$ has only finitely many isomorphic classes of 
simple modules and so it has only finitely many conformal weights 
$\{\la_1,..,\la_k\}$ and all of them are rational numbers. 
Therefore, if we can take an integer $n$ large enough ($n>\!\!>|\la_i-\la_j|$ for all $i,j$),  
then there is a bijective correspondence between 
composition factors of $W(n)$ as $A_n(V)$-modules and those of $W$ as $V$-modules for 
any module $W$. We then fix an integer $k$ satisfying $k>\max\{n+\la_i\mid i=1,...,k\}$. 
By the choice of $k$, for any indecomposable $V$-module $W$, 
$W$ has a nonzero homogeneous space $W_m=\{w\in W \mid (L(0)-m)^{\lambda}v=0 \mbox{ for some }\lambda\geq 0\}$ for some $m$ satisfying $k\leq m <k+1$. 
Take an integer $N\geq \max\{k-\la_i\mid i=1,...,k\}$. 
Then $W_m$ is an $A_{N}(V)$-module and there is a bijective correspondence 
between composition factors of $W_m$ as $A_N(V)$-modules and 
those of $W$ as $V$-modules. 
Therefore, from now on, we fix $k$ and $N$ and 
we will only consider $W_m$ satisfying $k\leq m<k+1$ as $A_N(V)$-modules 
and we denote it by $W[N]$. In other words, let $I_i$ be an ideal of $A_N(V)$ generated by 
$(\omega-\la_i-j)^s$ for a sufficiently large $s$ and an integer $j$ satisfying $k\leq \la_i+j<k+1$. 
Set $A_N(V)^{[k]}=A_N(V)/\sum I_i$ and we will consider only $A_N(V)^{[k]}$-modules. 
\end{rmk}

We now prove the finiteness of $\CF(W,U)$. 

\begin{prn} If $V$ is $C_2$-cofinite and $W$ and $U$ are f.g. $V$-modules then 
$\CF(W,U)$ at (2.5) contains only finitely many isomorphism classes.
\end{prn}

\noindent
[{\bf Proof}]\quad 
Suppose false, then since $\CF(W,U)$ is a direct set, it contains a strictly 
increasing infinite series $\CY_i\in \CI_{W,U}^{T^i}$
$$  (T^1,\CY_1)<(T^2,\CY_2)<(T^3,\CY_3)<\cdots. $$
By choosing a subsequence we may assume 
that the conformal weights $\wt(T^i)$ are all the same, say $t$, since 
$V$ has only finitely many of them. 
Let $N$ be a sufficiently large integer and we fix it in the proof.  
Set $\epsilon_{i,j}:T^i\rightarrow T^j$ for $i>j$ so that 
$\epsilon_{i,j}\cdot\CY_i=\CY_j$ and consider 
$$w^{\CY_i}_{(\wt(w)-1+\wt(U(0))-t,h)}:U(N) \rightarrow T^i(N)$$
for $w\in W$ and $h=0,1,...,K$, where $\CY_i(w,z)
=\sum_{h=0}^K\sum_{n\in \Q} w^{\CY_i}_{(n,h)}z^{-n-1}\log^hz$.
To simplify the notation, set 
$$o_*(\CY_{i}^{(h)})(w)=w^{\CY_i}_{(\wt(w)-1+\wt(U(0))-t,h)}.$$  
Then $\epsilon_{s,t}(o_*(\CY_{s}^{(h)})(w))=o_*(\CY_{t}^{(h)})(w)$. 
If we set 
$$  \widetilde{O}_{N,N}(W)
=<\sum_{i=0}^{\wt(v)+N}\binom{\wt(v)+N}{i}v^W_{-2N-2+i}w \mid v\in V, w\in W>, \eqno{(2.9)}$$
then by Borcherds' identity (2.2), we have 
$$o_*(\CY_{i}^{(h)})(w)=0$$
for any $w\in \widetilde{O}_{N,N}(W)$ and $i,h\in\N$. 
Therefore, there are epimorphisms 
$$\phi_i:\oplus_{h=0}^K\left((W/\widetilde{O}_{N,N}(W)\otimes U(N)\right) \rightarrow T^i(N)$$ 
given by $\phi_i(\oplus_{h=0}^K(w^h\otimes u^h))
=\sum_{h=0}^K o_*(\CY_{i}^{(h)})(w)u^h$ for any $i$. 

On the other hand, by the same arguments as in the proof 
in \cite{Miy}(Th 2.5) for the result \cite{GN}, 
if $B$ is a direct sum of homogeneous subspaces in $W$ satisfying  
$B+C_{2N+2}(W)=W$, then we can prove $B+\widetilde{O}_{N,N}(W)=W$ by (2.7). 
For, suppose false and take a homogeneous element 
$w\not\in B+\widetilde{O}_{N,N}(W)$ such that $\wt(w)$ is minimal.
Then since $B+C_{2N+2}(W)=W$, we obtain 
$w=b+\sum a^i_{-2N-2}c^i$ for some $b\in B$, 
$a^i\in V$, $c^i\in W$. Since $w$ is homogeneous and $B$ is a direct sum of 
homogeneous subspaces, we may assume $\wt(a^i_{-2N-2}c^i)=\wt(w)$ for 
all $i$. Furthermore, since 
$$a^i_{-2N-2}c^i+\sum_{j=1}^{\wt(a^i)+N}\binom{\wt(a^i)+N}{j}a^i_{-2N-2+j}c^i
\in \widetilde{O}_{N,N}(W)$$ 
for all $i$ and $\wt(a^i_{-2N-2+j}c^i)<\wt(w)$ for $j\geq 1$, 
we have $a^i_{-2N-2}c^i\in \widetilde{O}_{N,N}(W)+B$ by the minimality of $\wt(w)$. 
Therefore, we have $w\in \widetilde{O}_{N,N}(W)+B$, a contradiction and 
we have $B+\widetilde{O}_{N,N}(W)=W$. 
Since $V$ is $C_2$-cofinite, $W/C_{2N+2}(W)$ is of finite dimensional 
because of Proposition \ref{C2Pr} (iii) and so  
$W/\widetilde{O}_{N,N}(W)$ is also of finite dimensional. Hence a direct sum of 
$(K+1)$ copies of 
$(W/\widetilde{O}_{N,N}(W)\otimes U(N))$ is of 
finite dimensional and we have that 
$\dim(T^i(N))$ are bounded, which implies a contradiction. 
\prend

As a consequence, we have:

\begin{thm} 
If $V$ is a $C_2$-cofinite VOA, then for f.g. modules $W$ and $U$, 
there exists a tensor product $W\boxtimes U$ of $W$ and $U$ and it is 
also finitely generated. 
\end{thm}

We next explain about induced homomorphisms among tensor products. 
Let $(\CY^{B\boxtimes D},B\boxtimes D)$ be a tensor product of $B$ and $D$. 
For a homomorphism $\phi: A\rightarrow B$ of $V$-modules, a formal operator 
$\CY\in \CI_{A,D}^{B\boxtimes D}$ defined by 
$$\CY(a,z)d=\CY^{B\boxtimes D}(\phi(a),z)d$$ 
becomes an intertwining 
operator of type $\left({}_{A,D}^{B\boxtimes D}\right)$. Therefore, by the 
maximality of tensor products, there is a homomorphism 
$\phi\boxtimes {\rm id}_D: A\boxtimes D \rightarrow B\boxtimes D$ such that 
$$\phi\boxtimes {\rm id}_D\cdot(\CY(a,z)d)=\CY(\phi(a),z)d.$$  
We call $\phi\boxtimes {\rm id}_D$ an induced map of $\phi$.  
Similarly, we can define ${\rm id}_D\boxtimes \phi:D\boxtimes A\rightarrow D\boxtimes B$. 
From the definition, we easily have the right-flatness of modules 
for tensor products. 

\begin{prn}\label{RightFl}  Let $A$, $B$, $C$, $D$ be f.g. $V$-modules and 
assume that tensor products of f.g. modules are 
well-defined and are all finitely generated. Then if
$$A\xrightarrow{\phi} B\xrightarrow{\sigma} C \rightarrow 0$$ 
is exact, then so is 
$$ A\boxtimes D \xrightarrow{\phi\boxtimes {\rm id}_D} 
B\boxtimes D 
\xrightarrow{\sigma\boxtimes {\rm id}_D} C\boxtimes D \rightarrow 0.$$
\end{prn}


\pr
Clearly, $(\sigma\boxtimes {\rm id}_D)\cdot(\phi\boxtimes {\rm id}_D)
=(\sigma\cdot \phi)\boxtimes {\rm id}_D=0$. It is also clear that 
$\sigma\boxtimes {\rm id}_D$ is surjective and so we may view 
$(C\boxtimes D)'\subseteq (B\boxtimes D)'$. 
We may also assume $A\subseteq B$ and $C=B/A$. Consider a canonical bilinear pairing  
$$\langle g , \CY^{B\boxtimes D}(b,z)d\rangle \in \C\{z\}[\log z] $$
for $g\in (B\boxtimes D)', b\in B$ and $d\in D$. 
Clearly, if $g\in (C\boxtimes D)'$, then 
$$\langle g, \CY^{B\boxtimes D}(a,z)d\rangle =0$$
for any $a\in A$. 
On the other hand, if $g\in (B\boxtimes D)'$ satisfies 
$$\langle g,\CY^{B\boxtimes D}(a,z)d\rangle =0$$
for any $a\in A$, then $\langle g, \CY^{B\boxtimes D}(b,z)d\rangle $ is well 
defined for $b\in B/A=C$. Therefore, 
$$0 \rightarrow (C\boxtimes D)' \rightarrow (B\boxtimes D)' 
\rightarrow (A\boxtimes D)'$$
is exact, which implies the right flatness of $\boxtimes D$. 
\prend

\subsection{Projective covers}
Let us start with an explanation of projective modules and 
projective covers. 

\begin{dfn} A $V$-module $P$ is called projective 
if every $V$-epimorphism $f:W\rightarrow P$ is split. 
On the other hand, a $V$-module $Q$ is called injective  
if every $V$-monomorphism $g:Q\rightarrow W$ is split. 
\end{dfn}

Clearly, $P$ is projective if and only if $P'$ is injective. 
Different from an ordinary ring theory, 
$V$ is not necessary to be projective as a $V$-module. 
We first show the existence of projective covers.

\begin{prn}  Let $V$ be a $C_2$-cofinite VOA. Then we have: \\
(1) For any f.g. module $U$, there is a pair $(P,f)$ of 
a projective module $P$ and an epimorphism $f:P\rightarrow U$.  
If we assume that $\Ker f$ does not contain a direct summand of $P$, 
then $P$ is uniquely determined up to isomorphism. 
We call such a projective module $P$ a projective cover of $U$. \\
(2) Let $f:P\rightarrow U$ be a projective cover of $U$ and $g:E\rightarrow U$ an epimorphism. 
If $U$ is generated from $u$ as a $V$-module and $e\in E$ satisfies $g(e)=u$,  
then there are an integer $n$ and a homomorphism $h:P^{\oplus n}\rightarrow E$ such that 
$gh=\sum f^{(i)}$ and $e\in h(P^{\oplus n})$, 
where $P^{\oplus n}=\oplus_{i=0}^n P^{(i)}$, $P^{(i)}\cong P$, and 
$f^{(i)}:P^{(i)}\rightarrow U$ is a copy of $f:P\rightarrow U$. 
\end{prn}

\pr 
Since the uniqueness in (1) comes from ordinary ring theoretical arguments,
we will prove only the existence in (1) and then the statement (2).   
We may assume that $U$ is generated by one homogeneous element $u$. Set 
$$\CU=\left\{ (f,W,w) \mid \begin{array}{l}
\mbox{$W$ is a f.g. $V$-module}, w\in W, \wt(w)=\wt(u), Vw=W\cr
 f\in \Hom(W,U), f(w)=u \end{array} \right\},$$
where $Vw$ denotes a $V$-submodule generated from $w$. 
We introduce a partial order $(f',W',w')<(f,W,w)$ in $\CU$ 
if there is a surjective $V$-homomorphism $\phi:W\rightarrow W'$ 
such that $\phi(w)=w'$ and $f'\cdot\phi=f$. 

Since all $W$ are generated from one element $w$ with a fixed weight $\wt(u)$,  
$\dim W_n$ is bounded by Proposition \ref{C2Pr} (iii) and so 
the length of composition series of $W$ is bounded. 
Therefore, there is a maximal element $(\widetilde{f},\widetilde{W},\widetilde{w})$ in $\CU$. 
The next step is to show that $\widetilde{W}$ is projective. 
If there is an epimorphism $\phi:R\rightarrow \widetilde{W}$, then there is a homogeneous element 
$q\in R$ such that $\phi(q)=\widetilde{w}$. Since $V\widetilde{w}=\widetilde{W}$, 
$\phi:Vq\rightarrow \widetilde{W}=V\phi(q)$ is surjective.  
However, by the maximality of $\widetilde{W}$, 
it should be an isomorphism. Consequently, $\widetilde{W}$ is projective. \\
Proof of (2).  Suppose false and let $E$ be a counter example with a minimal 
length of composition series. Clearly, we have $E=Ve$. 
Then $(g,E,e)\in \CU$. If $E/{\rm rad}(E)$ is not simple, then there are homogeneous elements 
$e^i\in E$  ($i=1,...,k$) such that $e=\sum_{i=1}^k e^i$ and 
$Ve^i$ are all proper submodules of $E$. 
Then by the minimality of $E$, there are $P^{\oplus n_i}$ and $f^i:P^{\oplus n_i}\rightarrow Ve^i$ such that 
they satisfy the statements in the theorem. Then $\oplus_i (P^{\oplus n_i})$ and 
$\prod f^i$ satisfy the desired condition. Hence we may assume that 
$E/{\rm rad}(E)$ is simple. 
As we showed in the proof of (1), there is a $(\widetilde{f},\widetilde{W},\widetilde{w})$ in $\CU$ and 
$\mu:\widetilde{W}\rightarrow E$ such that $\widetilde{W}$ is a projective module and  
$\mu(\widetilde{w})=e$ and $g\mu=\widetilde{f}$. 
From the uniqueness of projective cover, $\widetilde{W}$ has $P$ as a direct summand 
such that $g\mu(P)=\widetilde{f}(P)=U$ and we may view $f=\widetilde{f}_{|P}$. 
In particular, $\mu(P)=E$ and $h:=\mu_{|P}$ is a desired homomorphism. 
\prend

A projective cover 
$$   0 \rightarrow {\rm rad}(P) \rightarrow P \xrightarrow{\rho} V \rightarrow 0  
\quad \mbox{ (exact)} $$
of a $V$-module $V$ will be called a "principal projective cover". 
From now on, we fix the notation $P$ to denote it.

If $V$ is $C_2$-cofinite and $M=\oplus_{i=0}^{\infty}M(i)$ is an $\N$-graded $V$-module, 
then as we mentioned in Remark 7, 
the composition series of $M$ as $V$-modules are corresponding to the composition series of 
$M[N]$ as $A_N(V)^{[k]}$-modules for sufficiently large integers $k$ and $N$.  
Furthermore, since there exists an 
$\N$-graded $V$-module $W$ satisfying $W[N]\cong R$ for 
any $A_N(V)^{[k]}$-module $R$, it is easy to see that $M$ is projective if and only if $M[N]$ 
is a projective $A_N(V)^{[k]}$-module. Furthermore, since $A_N(V)^{[k]}$ is a finite dimensional algebra, 
$M[N]$ is projective if and only if $M[N]$ is a direct sum of 
direct summands of a left $A_N(V)^{[k]}$-module $A_N(V)^{[k]}$.  
In particular, $M$ is a projective 
cover of a simple $V$-module $W$ if and only if $M[N]$ is an indecomposable 
direct summand of $A_N(V)^{[k]}$.  Therefore, we have the following:

\begin{prn} 
Let $V$ be a $C_2$-cofinite and $M$ an $\N$-graded $V$-module. 
Then the followings are equivalent. \\
(1)  $M$ is projective. \\
(2)  $M[N]$ is a direct sum of direct summands of a left $A_N(V)^{[k]}$-module $A_N(V)^{[k]}$. \\
\end{prn}

\section{Flatness}
\subsection{Key theorem}
In \S 2.4, we have explained a principal projective cover. 
As an application, we will construct every projective cover 
from it by using tensor products under the assumption of the semi-rigidity. 
Namely, we will prove:

\begin{thm}\label{KeyTh}
Assume that $V$ is a $C_2$-cofinite VOA of CFT type and all simple $V$-modules satisfy 
the semi-rigidity.  
Let $W$ be a f.g. $V$-module and 
$f:S \rightarrow W$ a projective cover of $W$.  
Then there is an epimorphism $g:W\boxtimes P \rightarrow S$ such that 
$fg$ coincides with ${\rm id}_W\boxtimes \rho:W\boxtimes P \rightarrow W\boxtimes V=W$, 
where $P$ is a principal projective cover. 
In particular, $S$ is a direct summand of $W\boxtimes P$. 
\end{thm}

\pr
If $W/{\rm rad}(W)=\oplus_i W^i$ and $S^i$ is 
a projective cover of $W^i$, then a projective cover of $W$ is a 
direct summand of $S^i$. Therefore, we may assume that $W$ is a simple $V$-module and 
$f:S\rightarrow W$ is its projective cover.  
Let $\mu:(W\boxtimes\widetilde{W})\boxtimes W \rightarrow W\boxtimes(\widetilde{W}\boxtimes W)$ 
be a canonical isomorphism. 
From the semi-rigidity, 
for $e_W:\widetilde{W}\boxtimes W\rightarrow V$, 
there is a submodule $R$ of $W\boxtimes \widetilde{W}$ with $R/{\rm rad}(R)\cong V$ 
such that $(id_W\boxtimes e_W)(\mu(R\boxtimes W))=W$. Since 
$S\boxtimes \widetilde{W}\rightarrow W\boxtimes \widetilde{W}$ is epimorphism, 
there is a submodule $Q$ of $S\boxtimes \widetilde{W}$ satisfying $Q/{\rm rad}(Q)\cong V$ 
such that $(f\boxtimes {\rm id}_{\widetilde{W}})(Q)=R$. Therefore, 
we have the following commutative diagram:
$$\begin{array}{rcccccc}
Q\boxtimes W&\subseteq &(S\boxtimes \widetilde{W})\boxtimes W &
\xrightarrow{\mu_S}& S\boxtimes (\widetilde{W}\boxtimes W) 
&\xrightarrow{{\rm id}_S\boxtimes e_W}& S\boxtimes V \cong S\cr
\quad \downarrow epi &   &\quad \downarrow epi     &&     \quad \downarrow epi    &&  \quad \downarrow epi  \cr
R\boxtimes W&\subseteq&(W\boxtimes \widetilde{W})\boxtimes W  &\xrightarrow{\mu}&
W\boxtimes (\widetilde{W}\boxtimes W) &\xrightarrow{{\rm id}_W\boxtimes e_W}
&W\boxtimes V\cong W
\end{array}$$
and $({\rm id}_W\boxtimes e_W)(\mu (R\boxtimes W))=W$. 
Since $S$ is a projective cover of $W$, $(id_S\boxtimes e_W)\mu_S:Q\boxtimes W\rightarrow S$ 
is surjective. 
By the choice of $Q$, $Q$ is a homomorphic image of a principal projective cover $P$ of $V$ and 
so we have a surjective homomorphism from $P\boxtimes W$ to $S$. 
Since $S$ is projective, $S$ is a direct summand of $P\boxtimes W$ as we desired.
\prend

\subsection{Flatness of modules for principal projective cover}
In this section, we will prove that 
for a principal projective cover, the tensor product with a 
f.g. module $F$ keeps an exactness.

\begin{prn}\label{PrinceFlat}
Let $V$ be a $C_2$-cofinite VOA of CFT type and $T$ is a 
f.g. $V$-module. Then 
$$0\rightarrow {\rm rad}(P)\boxtimes T 
\xrightarrow{\sigma\boxtimes {\rm id}_T} 
P\boxtimes T 
\xrightarrow{\rho\boxtimes {\rm id}_T} V
\boxtimes T\rightarrow 0  \eqno{(3.1)}$$
is exact, where $0\rightarrow {\rm rad}(P)\xrightarrow{\tau} P
\xrightarrow{\rho} V$ is a 
principal projective cover. 
\end{prn}

\pr  
Choose a homogeneous element 
$q\in P$ with $\rho(q)=\1$.  
Since
$$W=< v_nq \mid v\in V, n\in \Z >$$
is a submodule of $P$ by \cite{DM} and $\rho(W)=V$, 
we have $W=P$. Therefore, 
$${\rm rad}(P)=<v_nq\in P \mid v\in V, n\geq 0>. $$

Since the right flatness always holds, in order to prove the proposition, 
it is sufficient to show 
$$\sigma\boxtimes {\rm id}_T: {\rm rad}(P)\boxtimes T \rightarrow P\boxtimes T$$
is injective. More precisely, we will prove that ${\rm rad}(P)\boxtimes T$ is a 
submodule of $P\boxtimes T$ and 
$P\boxtimes T\cong T\oplus {\rm rad}(P)\boxtimes T$ as vector spaces. 
 Let $\CY^{{\rm rad}(P)\boxtimes T}=\sum_{j=0}^K\CY^{(j)}\log^jz$ 
be a tensor product intertwining operator. 
Then set $I=\CY^{(0)}$ and 
$I(w,z)=\sum_{r\in \Q} w^I_rz^{-r-1}$, 
where $w^I_r=w^{\CY^{{\rm rad}(P)\boxtimes T}}_{(r,0)}$ 
for $w\in {\rm rad}(P)$. We note that 
$I(w,z)$ is an $L(-1)$-nilpotent intertwining operator.   

Consider a vector space $R=q_{-1}T\oplus ({\rm rad}(P)\boxtimes T)$, where 
$q_{-1}T$ is a formal notation and $q_{-1}T\cong T$ as vector spaces. 
We will introduce a $V$-module structure on $R$ as follows:
View ${\rm rad}(P)\boxtimes T$ as a $V$-submodule and define an action of $v^R_n$ 
on $R=q_{-1}T\oplus {\rm rad}(P)\boxtimes T$ by 
$$v^R_n(q_{-1}t)
=q_{-1}(v_nt)+\sum_{i=0}^{\infty}\binom{n}{i}(v_iq)^I_{-1+n-i}t. \eqno{(3.2)}$$

By the direct calculation, we have Commutativity:
$$\begin{array}{l}
(u^R_mv^R_n-v^R_nu^R_m)(q_{-1}t) \cr
\mbox{}\quad=\displaystyle{q_{-1}((u_mv_n-v_nu_m)t)+
\sum_{i=0}^{\infty}\binom{m}{i}(u_iq)^I_{m-1-i}\ast v_nt
-\sum_{i=0}^{\infty}\binom{n}{i}(u_iq)^I_{m-1-i}v_nt}\cr
\mbox{}\qquad+\displaystyle{\sum_{i=0}^{\infty}\binom{n}{i}(v_iq)^I_{n-1-i}u_mt
-\sum_{i=0}^{\infty}\binom{n}{i}(v_iq)^I_{n-1-i}u_mt }\cr
\mbox{}\qquad+\displaystyle{\sum_{i=0}^{\infty}\sum_{j=0}^{\infty}
\binom{n}{i}\binom{m}{j}(u_jv_iq)^I_{m+n-1-i-j}t
-\sum_{i=0}^{\infty}\sum_{j=0}^{\infty}\binom{m}{j}\binom{n}{i}
(v_iu_jq)^I_{m+n-1-i-j}t }\cr
\mbox{}\quad=\displaystyle{q_{-1}([u_m,v_n]t)+\sum_{i=0}^{\infty}\sum_{j=0}^{\infty}
\binom{n}{i}\binom{m}{j}([u_i,v_j]q)^I_{m+n-1-i-j}t }\cr
\mbox{}\quad=\displaystyle{q_{-1}([u_m,v_n]t)+\sum_{i=0}^{\infty}\sum_{j=0}^{\infty}\sum_{k=0}^{\infty}
\binom{n}{i}\binom{m}{j}\binom{i}{k}((u_kv)_{i+j-k}q)^I_{m+n-1-i-j}t }\cr
\mbox{}\quad=\displaystyle{q_{-1}(\sum_{j=0}^{\infty}\binom{m}{j}(u_jv)^I_{m+n-j}t)
+\sum_{j=0}^{\infty}\sum_{i=0}^{\infty}\binom{m}{j}\binom{n+m-j}{i}((u_jv)_iq)^I_{m+n-1-j-i} }t\cr
\mbox{}\quad=\displaystyle{\sum_{j=0}^{\infty}\binom{m}{j}(u_jv)^R_{m+n-j}(q_{-1}t)}. 
\end{array}\eqno{(3.3)}$$
In particular, we have:
$$\begin{array}{rl}
\displaystyle{(\omega_0v)^R_n(q_{-1}t)}=&\displaystyle{q_{-1}((\omega_0v)_nt)+\sum_{i=0}^{\infty}
((\omega_0v)_iq)^I_{-1+n-i}t} \cr
=&\displaystyle{-nq_{-1}(v_{n-1}t)+\sum_{i=0}^{\infty}\binom{n}{i}(-iv_{i-1}q)^I_{-1+n-i}t} \cr
=&\displaystyle{-nv^R_{n-1}(q_{-1}t).}
\end{array}$$
For Associativity, we will not calculate it for every integer, but  
for only $m\in \Z_{\geq 0}$. 
By the direct calculation, we have:
$$\begin{array}{l}
\displaystyle{(v_nu)^R_m(q_{-1}t)-q_{-1}(v_nu)_mt}\cr
\mbox{}\quad=\displaystyle{\sum_{i=0}^{m}\binom{m}{i}((v_nu)_iq)^I_{m-1-i}t }\cr
\mbox{}\quad=\displaystyle{\sum_{i\in \N}\sum_{j\in \N}\binom{m}{i}\binom{n}{j}(-1)^j
(\{v_{n-j}u_{i+j}q-(-1)^nu_{n+i-j}v_jq\})^I_{m-1-1}t }\cr
\mbox{}\quad=\displaystyle{\sum_{i\in \N}\sum_{j\in \N}
\sum_{h\in \N}\binom{m}{i}\binom{n}{j}(-1)^{j+h}\binom{n-j}{h}
v_{n-j-h}(u_{i+j}q)^I_{m-i-1+h}t}\cr
\mbox{}\qquad-\displaystyle{\sum_{i\in \N}\sum_{j\in \N}\sum_{h\in\N}(-1)^{n+h}\binom{m}{i}\binom{n}{j}\binom{n-j}{h}
(u_{i+j}q)^I_{m-1+n-i-j-h}v_ht }\cr
\mbox{}\qquad+\displaystyle{\sum_{i\in \N}\sum_{j\in \N}\sum_{h\in \N}(-1)^{j+n+1+h}
\binom{m}{i}\binom{n}{j}\binom{n+i-j}{h}
u_{n+i-j-h}(v_jq)^I_{m-i-1+h}t}\cr
\mbox{}\qquad+\displaystyle{\sum_{i,j,h\in\N}(-1)^{i+h}
\binom{m}{i}\binom{n}{j}\binom{n+i-j}{h}
(v_jq)^I_{n+m-j-1-h}u_ht }\cr
\mbox{}\quad=\displaystyle{\sum_{i,j,h\in \N}\binom{m}{i}\binom{n}{j}(-1)^{j+h}\binom{n-j}{h}
\{A^{n-j-h}B^{i+j}\!-\!(-1)^{n-j+1}C^{i+j}D^h\} }\cr
\mbox{}\qquad+\displaystyle{\sum_{i,j,h\in\N}\binom{m}{i}(-1)^{j+n+1+h}\binom{n}{j}\binom{n+i-j}{h}
\{E^{n+i-j-h}F^j\!-\!(-1)^{n+i-j}G^jH^h\}}\cr
\mbox{}\quad=\displaystyle{(1+B)^m(A-1-B)^n-(1+C)^m(-1+D-C)^n}\cr
\mbox{}\qquad-\displaystyle{(1+E\!-\!1)^m(-E\!+\!1\!+\!F)^n+(1-1+H)^m(1-H+G)^n},   
\end{array}$$
where $A^sB^k$, $B^sD^k$, $E^sF^k$ and 
$G^sH^k$ denote $\displaystyle{v_{s}(u_{k}q)^I_{m+n-s-k-1}t}$, \\
$\displaystyle{(u_{s}q)^I_{m-1+n-s-k}v_kt}$, $\displaystyle{u_s(v_kq)^I_{m+n-k-s-1}t}$, and 
$\displaystyle{(v_sq)^I_{n+m-s-1-k}u_kt}$, respectively. 
On the other hand, for $n$-th normal product $\ast_n$, we have:
$$\begin{array}{l}
\displaystyle{(v \ast_n u)^R_m(q_{-1}t)-q_{-1}(v_n u)_mt }\cr
\mbox{}\quad
=\displaystyle{\sum_{h\in\N}\binom{n}{h}(-1)^h\{v_{n-h}u_{m+h}-(-1)^nu_{m+n-h}v_h\}(q_{-1}t)}\cr
\mbox{}\quad=\displaystyle{\sum_{h\in\N}\binom{n}{h}(-1)^hv_{n-h}\binom{h+m}{j}(u_jq)^I_{m+h-1-j}t}\cr
\mbox{}\qquad\displaystyle{+\sum_{h,j\in\N}\binom{n}{h}(-1)^h\binom{n-h}{j}
(v_jq)^I_{n-h-1-j}u_{m+h}t }\cr 
\mbox{}\qquad+\displaystyle{\sum_{h,j\in\N}\binom{n}{h}(-1)^{n+h+1}u_{m+n-h}
\binom{h}{j}(v_jq)^I_{h-1-j}t }\cr
\mbox{}\qquad+\displaystyle{\sum_{h,j\in\N}\binom{n}{h}(-1)^{n+h+1}
\binom{m+n-h}{j}(u_jq)^I_{n+m-h-j}v_ht }\cr
\mbox{}\quad=\displaystyle{\sum_{h,j\in\N}\binom{n}{h}(-1)^h\binom{m+h}{j}A^{n-h}B^j
+\sum_{h,j\in\N}\binom{n}{h}(-1)^{h+m+1}\binom{n+m-h}{j}G^jH^h }\cr
\mbox{}\qquad+\displaystyle{\sum_{h,j\in\N}\binom{n}{h}(-1)^{n+h+1}\binom{h}{j}E^{n+m-h}F^j
+\sum_{h,j\in\N}\binom{n}{h}(-1)^{h}\binom{n-h}{j}C^jD^{m+h}}\cr
\mbox{}\quad=\displaystyle{(1+B)^m(A-1-B)^n+H^m(1+G-H)^n-E^m(-E+1+F)^n}\cr
\mbox{}\qquad\displaystyle{-(1+C)^m(-1-C+D)^n}.
\end{array}$$ 
Therefore, we have Associativity 
$(v\ast_nu)^R_m r=(v_nu)^R_m r$ for $r\in R$ and $m\geq 0$. 

Define 
$$v^R(z)=\sum_{n\in \Z} v^R_n z^{-n-1}\qquad v^R_n \in \End(R)$$
for $v\in V$. Since the weights of $R$ are bounded below, 
$v^R(z)$ satisfies the lower truncation property. 
Furthermore they are mutually commutative because of (3.3) and so  
they generate a local system $\widetilde{V}$ in $\End(R)[[z,z^{-1}]]$ 
by using normal products $\ast_n$. 
Set $\overline{V}=\{v^R(z)\mid v\in V\}\subseteq \widetilde{V}$. 
Since we have  
$$\Res_{x}\{(x-z)^mv^R(x)u^R(z)-(-z+x)^mu^R(z)v^R(x)\}=(v_mu)^R(z)$$
for $m\geq 0$, $\overline{V}$ is closed under $m$-th product for $m\geq 0$. 
Furthermore, since \\
$[\omega^R_0,v^R_n]=-nv^R_{n-1}$ for Virasoro element $\omega$ of $V$, 
$\omega^R(z)$ is a Virasoro element of $\widetilde{V}$. 
Hence the grading of $\widetilde{V}$ is the same as on $\overline{V}$. 
We will prove $\overline{V}=\widetilde{V}$. 
Suppose false, then $\overline{V}$ is not closed in $\widetilde{V}$ by $n$-th product 
for some $n<0$. In other words, there are $v^R(z),u^R(z)\in \overline{V}$ 
such that $v^R(z)_n\ast u^R(z)\not\in \overline{V}$ for some $n$. 
We take $v^R(z)_n\ast u^R(z)\not\in \overline{V}$ such that 
$\wt(v^R(z)_n\ast u^R(z))$ is minimal. 
Set $\alpha(z)=v^R(z)_n\ast u^R(z)-(v_nu)^R(z)$ and $N=\wt(\alpha(z))$. 
We note $L(k)(\alpha(z))=0$ for $k\geq 1$ because of the 
minimality of the weight of $\alpha(z)$.   
Since $\alpha(z)_m=0$ on $R$ for $m\geq 0$ as we showed, we have: 
$$0=[\omega^R(z)_k,\alpha(z)_m]
=(\omega_0\alpha(z))_{k+m}+kN\alpha(z)_{k+m-1}=(-k-m+kN)\alpha(z)_{k+m-1}$$
for any $k\in \Z$ and so $\alpha(z)_k=0$ for any $k$, 
which contradict to $\alpha(z)\not=0$. 
Therefore, $\overline{V}$ is closed by $n$-th product for any $n\in \Z$ and 
$R$ becomes a $V$-module.

We next define an intertwining operator 
of $P$ from $T$ to $R$ by using $I$. 
We first note that since 
$$\begin{array}{l}
\displaystyle{v^R_n(L(-1)q)^I_{0}-(L(-1)q)^I_{0}v^R_n} \hspace{9cm} \cr
\mbox{}\hspace{4cm}=\displaystyle{\sum_{i=0}^{\infty}\binom{n}{i}(v^P_iL(-1)q)^I_{n-i} }\cr
\mbox{}\hspace{4cm}=\displaystyle{\sum_{i=0}^{\infty}\binom{n}{i}\{(L(-1)v^P_iq)^I_{n-i}+(iv^P_{i-1}q)^I_{n-i}\} }\cr
\mbox{}\hspace{4cm}=\displaystyle{\sum_{i=0}^{\infty}\binom{n}{i}\{(-n+i)(v^P_iq)^I_{n-i-1}+(iv^P_{i-1}q)^I_{n-i}\}}\cr
\mbox{}\hspace{4cm}=\displaystyle{-n\sum_{i=0}^{\infty}\binom{n-1}{i}(v^P_iq)^I_{n-i-1}
+\sum_{i=0}^{\infty}n\binom{n-1}{i-1}(v^P_{i-1}q)^I_{n-i}}\cr
\mbox{}\hspace{4cm}=0,
\end{array}$$
$(L(-1)q)^I_{0}$ commutes with the actions of $V$. 
We define $q(z):T\rightarrow R[[z,z^{-1}]]$ by 
$$  q(z)t=\left\{\int (\CY^{{\rm rad}(P)\boxtimes T}(L(-1)q,z)-(L(-1)q)_{0}z^{-1})dz\right\}t
+q_{-1}t. \eqno{(3.4)}$$ 
Let us show that $q(z)$ satisfies Commutativity with all actions 
$a(z)$ of $a\in V$. 
By direct calculation, we have: 
$$\begin{array}{l}
[a(z), q(x)]-\dsp{\sum_{n\in \Z}\sum_{i=0}^{\infty}\binom{n}{i}(a_iq)^I_{n-1-i}z^{-n-1}} \cr
\mbox{}\qquad=\dsp{\sum_{m\not=0,n\in \Z}[a_nz^{-n-1}, (L(-1)q)^I_m
\frac{-1}{m}x^{-m}] }\cr
\mbox{}\qquad=\dsp{\sum_{m\not=0,n\in \Z}\sum_{i=0}^{\infty}\binom{n}{i}(a_iL(-1)q)^I_{m+n-i}
\frac{-1}{m}x^{-m}z^{-n-1}} \cr
\mbox{}\qquad=\dsp{\sum_{m\not=0,n\in \Z}\sum_{i=0}^{\infty}\binom{n}{i}\{(L(-1)a_iq)^I_{m+n-i}
+i(a_{i-1}q)^I_{m+n-i}\}
\frac{-1}{m}x^{-m}z^{-n-1}}\cr
\mbox{}\qquad=\dsp{\sum_{m\not=0,n\in \Z}\sum_{i=0}^{\infty}\binom{n}{i}\{(-m-n+i)(a_iq)^I_{m+n-i-1}
+i(a_{i-1}q)^I_{m+n-i}\}\frac{-1}{m}x^{-m}z^{-n-1}} \cr
\mbox{}\qquad=\dsp{\sum_{m\not=0,n\in \Z}\sum_{i=0}^{\infty}\binom{n}{i}(-m)(a_iq)^I_{m+n-i-1}
\frac{-1}{m}x^{-m}z^{-n-1}}\cr
\mbox{}\qquad=\dsp{\sum_{m\not=0,n\in \Z}
\sum_{i=0}^{\infty}\binom{n}{i}(a_iq)^I_{m+n-i-1}x^{-m}z^{-n-1}},
\end{array}$$
and so 
$$[a(z),q(x)]=\sum_{m,n\in \Z}\sum_{i\in \N}
\binom{n}{i}(a_iq)^{I}_{(m+n-i-1)}x^{-m}z^{-n-1}.$$
If we take $N$ so that $a_iq=0$ for $i\geq N$, then 
$$\begin{array}{rl}
[a(z), q(x)]=&\dsp{\sum_{m,n\in \Z}\sum_{i\in \N}\binom{n}{i}(a_iq)^I_{n+m-1-i}x^{-m}z^{-n-1}} \cr
=&\dsp{\sum_{i=0}^N\sum_{r\in \Z}(a_iq)^I_{r-i-1}\sum_{n\in \Z}\binom{n}{i}x^{n-r}z^{-n-1}} 
\end{array}$$
and so we have locality:
$$(x-z)^{N+1}[a(z),q(x)]=0.$$ 
We then extend it by  
$$J^{(0)}(v_nq,z)=\Res_{x}\{(x-z)^nv(x)q(z)-(-z+x)^nq(z)v(x)\}$$
for $v\in V, n\in \Z$, then 
$J^{(0)}$ is an $L(-1)$-nilpotent intertwining operator and so 
$$J(w,z)=\sum_{i=0}^K\frac{1}{i!}(zL(-1)-z\frac{d}{dz})^iJ^{(0)}(w,z)\log^iz $$
is a surjective intertwining operator of $P$ from $T$ to $R$. On the other hand, 
since the right flatness holds, in other words, since there is an exact sequence 
${\rm rad}(P)\boxtimes T \rightarrow P\boxtimes T \rightarrow V\boxtimes T\rightarrow 0$, we obtain 
from the length of compositions series that  
$${\rm rad}(P)\boxtimes T \xrightarrow{\sigma\boxtimes {\rm id}_T} 
P\boxtimes T$$ is injective. 
\prend

\subsection{Main theorem}
We now start the proof of the main theorem. \\

\begin{thm}\label{MainTh} Let $V$ be a $C_2$-cofinite vertex operator algebra of CFT type and 
$A$, $B$, $C$, $D$ f.g. $V$-modules. 
Assume that all simple $V$-modules satisfy the semi-rigidity. If 
$$0 \rightarrow A \xrightarrow{\tau} B 
\xrightarrow{\sigma} C\rightarrow 0$$ 
is an exact sequence of $V$-modules, then so is 
$$0\rightarrow A\boxtimes D \xrightarrow{\tau\boxtimes {\rm id}_D} 
B\boxtimes D 
\xrightarrow{\sigma\boxtimes {\rm id}_D}C\boxtimes D \rightarrow 0. \eqno{(3.5)}$$
\end{thm}

\pr
Since we may assume that ${\rm Ker}~\sigma$ does not contains 
a direct summand of $B$, we have a short exact sequence 
$$0\rightarrow J\rightarrow P\boxtimes C \rightarrow B \rightarrow 0$$
from Theorem \ref{KeyTh}.  
Using this and (3.1), we have the following commutative exact diagram. 
$$\begin{array}{ccccccccc}
 &&     0      &\rightarrow&      0      &&          &&  \cr
      &  &\dto&  &\dto& && & \cr
0 &\rightarrow&     J      &\rightarrow&      J      &\rightarrow&    0       &&  \cr
      \dto&  &\dto&  &\dto& &\dto& & \cr
0 &\rightarrow &{\rm rad}(P)\boxtimes C &\rightarrow& P\boxtimes C  &\rightarrow& V\boxtimes C=C &\rightarrow&  0 \cr
      \dto&  &\dto&  &\dto& &\dto& &\dto\cr
0 &\rightarrow&     A      &\rightarrow&      B      &\rightarrow&    C       &\rightarrow&  0 \cr
      &  &\dto&  &\dto& &\dto& &\cr
 &&     0      &\rightarrow&      0      &\rightarrow&    0       &&   
\end{array}$$ 
Taking a tensor product of $C\boxtimes D$ with a principal projective cover, 
we have an exact sequence 
$$  0\rightarrow {\rm rad}(P)\boxtimes (C\boxtimes D) \rightarrow  
P\boxtimes (C\boxtimes D) \rightarrow C\boxtimes D\rightarrow 0$$
by $(3.1)$. Using isomorphisms  
$A\boxtimes (B\boxtimes C)\cong (A\boxtimes B)\boxtimes C$
and from the length of composition series, we have 
an exact sequence   
$$  0\rightarrow ({\rm rad}(P)\boxtimes C)\boxtimes D \rightarrow  
(P\boxtimes C)\boxtimes D \rightarrow C\boxtimes D\rightarrow 0.  $$
Therefore, using the right flatness of $\boxtimes D$, we have 
the following commutative exact diagram:
$$\begin{array}{ccccccccc}
0 &\rightarrow&     J\boxtimes D    &\rightarrow&      J\boxtimes D      &\rightarrow&    0       &&  \cr
    \dto  &  &\dto&  &\dto& &\dto& & \cr
0 &\rightarrow &{\rm rad}(P)\boxtimes C\boxtimes D &\rightarrow& P\boxtimes C\boxtimes D  &\rightarrow
  & C\boxtimes D &\rightarrow&  0 \cr
  &  &\dto&  &\dto& &\dto& &\dto\cr
  & &     A\boxtimes D      &\rightarrow&      B\boxtimes D      &\rightarrow
  &    C\boxtimes D       &\rightarrow&  0 \cr
  &  &\dto&  &\dto& &\dto& &\cr
  & &     0    &\rightarrow&      0      &\rightarrow&    0       &&  
\end{array}$$ 
which implies that $\sigma\boxtimes {\rm id}_D:A\boxtimes D\rightarrow B\boxtimes D$ 
is injective. 

This completes the proof of the main theorem. 
\prend


\begin{thebibliography}{99}
\bibitem[1]{B}
R.E.~Borcherds, {\it Vertex algebras, Kac-Moody algebras, and the 
Monster}, \\
{\it Proc.~Natl.~Acad.~Sci. USA} {\bf 83} (1986), 3068--3071.

\bibitem[2]{Bu}
G.~Buhl, {\it A spanning set for VOA modules}, J. Algebra {\bf 254}, 
(2002), 125-151.

\bibitem[3]{DLM} 
C.~Dong, H.~Li, and G.~Mason, 
{\it Twisted representations of vertex operator algebras and 
associative algebras},
Internat. Math. Res. Notices {\bf 8} (1998), 389--397.

\bibitem[4]{DM}
C.~Dong and G.~Mason,  {\it On quantum Galois theory}, Duke Math. J. 
{\bf 86} (1997), no.2, 305--321.

\bibitem[5]{F}
M.A.~Flohr, On modular invariant partition functions of conformal 
field theories with logarithmic operators, 
{\it Int. J. Mod. Phys.} A. (1995)

\bibitem[6]{FHL} 
I.~Frenkel, Y.-Z.~Huang and J.~Lepowsky, On axiomatic approaches to
vertex operator algebras and modules, \textit{Mem.~Amer.~Math.~Soc.} 
{\bf104} (1993).

\bibitem[7]{FLM}
I.~Frenkel, J.~Lepowsky, and A.~Meurman, "Vertex Operator Algebras 
and the Monster", Pure and Applied Math., Vol. 134, 
Academic Press, 1988.


\bibitem[8]{G} V.~Gurarie, {\it Logarithmic operators in conformal field theory}, 
Nuclear Phys. B {\bf 410} (1993), 535-549. 



\bibitem[9]{GN} M.~Gaberdiel and A.~Neitzke, Rationality, 
{\it quasirationality, and finite $W$-algebra}, DAMTP-200-111. 

\bibitem[10]{H1}
Y.-Z.~Huang, {\it Differential equations, duality and modular invariance},  
Commun. Contemp. Math. 7 (2005), no. 5, 649-706. 

\bibitem[11]{H2}
Y.-Z.~Huang, {\it Vertex operator algebras and the Verlinde conjecture},  
Commun. Contemp. Math. 10 (2008), no. 1, 103-154.   


\bibitem[12]{HL} Y.-Z.~Huang and J.~Lepowsky, {\it A theory of tensor products for 
module categories for a vertex algebra, I}, 
Secta Mathematica, New Series Vol. {\bf 1} (1995), 699-756.

\bibitem[13]{HLZ} Y.-Z.~Huang, J.~Lepowsky and L.~Zhang {\it A logarithmic generalization 
of tensor product theory for moduls for a vertex operator algebra}, 
Internat. J. Math. 17 (2006), no. 8, 975-1012. 

\bibitem[14]{Miy} M.~Miyamoto, 
{\it Modular invariance of vertex operator algebra satisfying 
$C_2$-cofiniteness}, Duke Math. J. 122 (2004), no. 1, 51-91. 

\bibitem[15]{Mil} A.~Milas, {\it Weak modules and logarithmic 
intertwining operators for vertex operator algebras}, Contemp. 
Math.,{\bf 297}, (2002) 201-225.

\bibitem[16]{TK} A.~Tsuchiya and Y.~Kanie, 
{\it Vertex operators in conformal 
field theory on ${\mathbb P}^1$ and monodromy representations 
of braid group}, in: Conformal Field Theory and Solvable Lattice 
Models, Advanced Studies in Pure Math. 16, Academic Press (1988), 297-372.

\bibitem[17]{Z} Y.~Zhu, {\it Modular invariance of characters of 
vertex operator algebras},  J. Amer. Math. Soc., {\bf 9} (1996), 
237--302. 
\end{thebibliography}
\end{document}